\documentclass{elsart}
\usepackage{amsfonts}
\usepackage{amsmath}

\input{tcilatex}

\begin{document}

\begin{frontmatter}%

\title{On the structure of the adjacency matrix of the line digraph of a regular digraph}%

\author{Simone Severini}%

\address{Department of Computer Science, University of Bristol, Merchant Venturers' Building, Woodland road, Bristol BS8 1UB, United Kingdom}%

\begin{abstract}
We show that the adjacency matrix $M$ of the line digraph of a $d$-regular
digraph $D$ on $n$ vertices can be written as $M=AB$, where the matrix $A$
is the Kronecker product of the all-ones matrix of dimension $d$ with the
identity matrix of dimension $n$ and the matrix $B$ is the direct sum of the
adjacency matrices of the factors in a dicycle factorization of $D$. 
\end{abstract}%


\end{frontmatter}%

\section*{Introduction}

Line digraphs of regular digraphs and their generalizations are important in
the design of point-to-point interconnection networks for parallel computers
and distributed systems. For instance, de Bruijn digraphs and
Reddy-Pradhan-Kuhl digraphs, which are important topologies for
interconnection networks, are all examples of line digraphs of regular
digraphs (see, \emph{e.g.}, \cite{BG00},\cite{Fe99} and \cite{H97}). In this
note, we describe a special regularity property of the adjacency matrix of
the line digraph of a regular digraph. Before stating formally our main
result, we recall the necessary graph-theoretic terminology.

A (\emph{finite}) \emph{directed graph}, for short \emph{digraph}, consists
of a non-empty finite set of elements called \emph{vertices} and a (possibly
empty)\ finite set of ordered pairs of vertices called \emph{arcs}. The
digraphs considered here are without multiple arcs. We denote by $D=\left(
V,A\right) $ a digraph with vertex-set $V(D)$ and arc-set $A\left( D\right) $%
. A \emph{labeling} of the vertices of a digraph $D$ is a function $%
l:V\left( D\right) \longrightarrow L$, where $L$ is a set of labels. Chosen
a bijective labeling, the \emph{adjacency matrix} of a digraph $D$ with $n$
vertices, denoted by $M\left( D\right) $, is the $n\times n$ $\left(
0,1\right) $-matrix with $ij$-th element defined by $M_{i,j}\left( D\right)
=1$ if $\left( v_{i},v_{j}\right) \in A\left( D\right) $ and $M_{i,j}\left(
D\right) =0$, otherwise. For any vertex $v_{i}\in V(D)$ of a digraph $D$,
let $d_{D}^{-}\left( v_{i}\right) :=|\left\{ v_{j}:\left( v_{j},v_{i}\right)
\in A\left( D\right) \right\} |$ and $d_{D}^{+}\left( v_{i}\right)
:=|\left\{ v_{j}:\left( v_{i},v_{j}\right) \in A\left( D\right) \right\} |$.
A digraph $D$ is said to be $d$\emph{-regular} if, for every vertex $%
v_{i}\in V(D)$, $d_{D}^{-}\left( v_{i}\right) =d_{D}^{+}\left( v_{i}\right)
=d$. A digraph $H$ is a \emph{subdigraph} of a digraph $D$ if $V\left(
H\right) \subseteq V\left( D\right) $ and $A\left( H\right) \subseteq
A\left( D\right) $. A subdigraph $H$ of a digraph $D$ is said to be a \emph{%
spanning subdigraph} of $D$, or equivalently, a \emph{factor} of $D$, if $%
V\left( H\right) =V\left( D\right) $. A \emph{decomposition} of a digraph $D$
is a set $\left\{ H_{1},H_{2},...,H_{k}\right\} $ of subdigraphs of $D$
whose arc-sets are exactly the classes of a partition of $A\left( D\right) $%
. A \emph{factorization} of a digraph $D$, if there exists one, is a
decomposition of $D$ into factors. A \emph{dicycle factor} $H$ of a digraph $%
D$ is a spanning subdigraph of $D$ such that $M\left( H\right) $ is a
permutation matrix. The \emph{disjoint union} of digraphs $%
D_{1},D_{2},...,D_{k}$, is the digraph with vertex-set $%
\biguplus_{i=1}^{k}V(D_{i})$, and arc-set $\biguplus_{i=1}^{k}A(D_{i})$.
Then a dicycle factor $H$ of a digraph $D$ is a spanning subdigraph of $D$
and it is the disjoint union of dicycles. A \emph{dicycle factorization} is
a factorization into dicycle factors. The \emph{line digraph} of a digraph $%
D $, denoted by $\overrightarrow{L}D$, is defined as follows: the vertex-set
of $\overrightarrow{L}D$ is $A\left( D\right) $; for $%
v_{h},v_{i},v_{j},v_{k}\in V\left( D\right) $, $(\left( v_{h},v_{i}\right)
,\left( v_{j},v_{k}\right) )\in A(\overrightarrow{L}D)$ if and only if $%
v_{i}=v_{j}$. Kronecker product and direct sum of matrices $M$ and $N$ are
respectively denoted by $M\otimes N$ and $M\oplus N$. The identity matrix
and the all-ones matrix of size $n$ are respectively denoted by $I_{n}$ and $%
J_{n}$. In the next section, we prove the following theorem:

\bigskip

\noindent \textbf{Theorem} \emph{Let }$D$\emph{\ be a }$d$\emph{-regular
digraph on }$n$\emph{\ vertices and let }$\left\{
H_{1},H_{2},...,H_{d}\right\} $\emph{\ be a dicycle factorization of }$D$%
\emph{. Then there is a labeling of }$V(\overrightarrow{L}D)$\emph{\ such
that }%
\begin{equation*}
M(\overrightarrow{L}D)=(J_{d}\otimes
I_{n})\dbigoplus\limits_{i=1}^{d}M(H_{i}). 
\end{equation*}

\section{Proof of the theorem}

The proof of the theorem is based on two simple observations and a result
proved by Hasunuma and Shibata \cite{HS96} (see also Kawai \emph{et al.} 
\cite{KFS01}).

\begin{lemma}
Let $D$ be a $d$-regular digraph. Then $D$ has a dicycle factorization. In
particular, if $\left\{ H_{1},H_{2},...,H_{d}\right\} $ is a dicycle
factorization of $D$ then $M(H_{1}),M(H_{2}),...,M(H_{d})$ are permutation
matrices such that%
\begin{equation*}
M(D)=\dsum\limits_{i=1}^{d}M(H_{i}). 
\end{equation*}
\end{lemma}

Two digraphs $D$ and $D^{\prime }$ are said to be \emph{isomorphic} if there
is a permutation matrix $P$ such that $P\cdot M(D)\cdot P^{-1}=M(D^{\prime
}) $. If $D$ and $D^{\prime }$ are isomorphic we then write $D\cong
D^{\prime }$. An $n$\emph{-dicycle}, denoted by $\overrightarrow{C}_{n}$, is
a digraph with vertex-set $\{v_{1},v_{2},...,v_{n}\}$ and arc-set $%
\{(v_{1},v_{2}),...,(v_{n-1},v_{n}),(v_{n},v_{1})\}$. A $d$-\emph{spiked }$n$%
\emph{-dicycle} is the digraph obtained from $\overrightarrow{C}_{n}$ as
follows: for every vertex $v_{i}\in V(\overrightarrow{C}_{n})$, we add $d$
new vertices $w_{1},w_{2},...,w_{d}$; we connect $v_{i}\in (\overrightarrow{C%
}_{n})$ to the vertices $w_{1},w_{2},...,w_{d}$, obtaining the arcs $\left(
v_{i},w_{1}\right) ,\left( v_{i},w_{2}\right) ,...,\left( v_{i},w_{d}\right) 
$.

\begin{lemma}
\label{cle}Let $D$ be a $d$-spiked $n$-dicycle. Then $D\cong \overrightarrow{%
L}D$.
\end{lemma}

Let $D$ be a digraph and let $H$ be a spanning subdigraph of $D$. The \emph{%
growth} of $D$\ derived by $H$ is the digraph denoted by $\Upsilon
_{D}\left( H\right) $ and defined as follows: for every pair of vertices $%
v_{i},v_{j}\in V(D)$, if $(v_{i},v_{j})\in A(H)$ then $(v_{i},v_{j})\in
A(\Upsilon _{D}\left( H\right) )$; for every vertex $v_{i}\in V(D)$, we add
new vertices $w_{1},w_{2},...,w_{l}$, where $l=d_{D}^{+}\left( v_{i}\right)
-d_{H}^{+}\left( v_{i}\right) $; we connect $v_{i}\in V(D)$ to the vertices $%
w_{1},w_{2},...,w_{l}$, obtaining the arcs $\left( v_{i},w_{1}\right)
,\left( v_{i},w_{2}\right) ,...,\left( v_{i},w_{l}\right) $.

\begin{lemma}[\protect\cite{HS96}]
\label{sht}If $\left\{ H_{1},H_{2},...,H_{k}\right\} $ is a decomposition of
a digraph $D$ then 
\begin{equation*}
\{\overrightarrow{L}\Upsilon _{D}(H_{1}),\overrightarrow{L}\Upsilon
_{D}(H_{2}),...,\overrightarrow{L}\Upsilon _{D}\left( H_{k}\right) \} 
\end{equation*}%
is a decomposition of a digraph $D^{\prime }\cong \overrightarrow{L}D$.
\end{lemma}

\begin{proof}[Proof of the theorem]
Let $D$ be a $d$-regular digraph on $n$ vertices $v_{1},v_{2},...,v_{n}$.
Let $\left\{ H_{1},H_{2},...,H_{d}\right\} $ be a dicycle factorization of $%
D $. The vertices of $H_{j}\in \left\{ H_{1},H_{2},...,H_{d}\right\} $ are
denoted as $(H_{j},v_{1}),(H_{j},v_{2}),...,(H_{j},v_{n})$. Let us construct 
$\Upsilon _{D}\left( H_{j}\right) $. For every vertex $(H_{j},v_{i})\in
V(H_{j})$, we add $d-1$ new vertices to $H_{j}$. We label these new vertices
by pairs of the form $(H_{l},v_{m})$, for all $l\neq j$ and $v_{m}$ such
that $(v_{i},v_{m})\in A(H_{l})$. In addition, $%
((H_{j},v_{i}),(H_{l},v_{m}))\in A(\Upsilon _{D}\left( H_{j}\right) )$. The
digraph $\Upsilon _{D}\left( H_{j}\right) $ has $n\cdot d$ vertices. If we
label the row number $(j-1)n+i$ of $M(\Upsilon _{D}\left( H_{j}\right) )$ by
the vertex $(H_{j},v_{i})$, the adjacency matrix of $\Upsilon _{D}\left(
H_{j}\right) $ is the $(d\cdot n)\times (d\cdot n)$ block-matrix%
\begin{equation*}
M(\Upsilon _{D}\left( H_{j}\right) )=\left( 
\begin{array}{c}
\mathbf{0} \\ 
X_{j} \\ 
\mathbf{0}%
\end{array}%
\right) , 
\end{equation*}%
where%
\begin{equation*}
X_{j}=\left( 
\begin{array}{ccccccc}
M\left( H_{1}\right) & M\left( H_{2}\right) & \cdots & M\left( H_{j}\right)
& \cdots & M\left( H_{d-1}\right) & M\left( H_{d}\right)%
\end{array}%
\right) . 
\end{equation*}%
Notice that $M\left( H_{j}\right) $ is the $jj$-th block of $M(\Upsilon
_{D}\left( H_{j}\right) )$. Thus, we have%
\begin{eqnarray*}
N &=&\dsum\limits_{i=j}^{d}M(\Upsilon _{D}\left( H_{i}\right) )=\left( 
\begin{array}{cccc}
M\left( H_{1}\right) & M\left( H_{2}\right) & \cdots & M\left( H_{d}\right)
\\ 
M\left( H_{1}\right) & M\left( H_{2}\right) & \cdots & M\left( H_{d}\right)
\\ 
\vdots & \vdots & \ddots & \vdots \\ 
M\left( H_{1}\right) & M\left( H_{2}\right) & \cdots & M\left( H_{d}\right)%
\end{array}%
\right) \\
&=&(J_{d}\otimes I_{n})\dbigoplus\limits_{i=j}^{d}M(H_{j}).
\end{eqnarray*}%
Observe that, for every $1\leq j\leq d$, $\Upsilon _{D}\left( H_{j}\right) $
is the disjoint union of the $d$-spiked cycles corresponding to the orbits
of the permutation associated to $H_{j}$. It follows from Lemma \ref{cle}
that, for every $1\leq j\leq d$, 
\begin{equation*}
\Upsilon _{D}(H_{j})\cong \overrightarrow{L}\Upsilon _{D}(H_{j}). 
\end{equation*}%
Then, for the chosen labeling,%
\begin{equation*}
M(\Upsilon _{D}\left( H_{j}\right) )=M(\overrightarrow{L}\Upsilon
_{D}(H_{j})) 
\end{equation*}%
and%
\begin{equation*}
N=\dsum\limits_{j=1}^{d}M(\Upsilon _{D}\left( H_{j}\right)
)=\dsum\limits_{j=1}^{d}M(\overrightarrow{L}\Upsilon _{D}(H_{j})). 
\end{equation*}%
Now, by Lemma \ref{sht}, $N=M(\overrightarrow{L}D)$.
\end{proof}

\bigskip

\noindent \textbf{Remark }The graph operation transforming a digraph $D$ in
its line digraph can be naturally iterated: $\overrightarrow{L}^{k}D:=%
\overrightarrow{L}\overrightarrow{L}^{k-1}D$. Let $\Sigma $ be an alphabet
of cardinality $d$ and let $\Sigma ^{k}$ be the set of all the words of
length $k$ over $\Sigma $. The $d$\emph{-ary }$k$\emph{-dimensional de Bruijn%
} \emph{digraph}, denoted by $B(d,k)$, is defined as follows: the vertex-set
of $B(d,k)$ is $V(B(d,k))=\Sigma ^{k}$; for every pair of vertices $%
v_{i},v_{j}$, we have $(v_{i},v_{j})\in A(B(d,k))$ if and only if the last $%
k-1$ letters of $v_{i}$ are the same as the first $k-1$ letters of $v_{j}$.
Let $K_{d}^{+}$ be the complete digraph on $d$ vertices with a loop at each
vertex. Fiol, Yebra and Alegre \cite{FYA84} proved that $B\left( d,k\right)
\cong \overrightarrow{L}^{k-1}K_{d}^{+}$. This result, together with the
theorem, gives 
\begin{equation*}
M(B(d,2))\cong (J_{d}\otimes I_{d})\dbigoplus\limits_{i=1}^{d}M(H_{i}), 
\end{equation*}%
where $\left\{ H_{1},H_{2},...,H_{d}\right\} $ is any dicycle factorization
of $K_{d}^{+}$.

\bigskip

\noindent \emph{Acknowledgments }The author would like to thank the
anonymous referees for helpful comments that improved this paper.


\begin{thebibliography}{9}
\bibitem{BG00} J. Bang-Jensen and G. Gutin, \emph{Digraphs. Theory,
algorithms and applications}, Springer Monographs in Mathematics,
Springer-Verlag, London, 2001.

\bibitem{Fe99} D. Ferrero, Introduction to interconnection network models, 
\emph{Publ.Mat. Urug.}, 99/25 (1999).

\bibitem{FYA84} M. A. Fiol, J. L. A. Yebra and I. Alegre, Line digraph
iterations and the $\left( d,k\right) $ digraph problem, \emph{IEEE Trans.
Comput.} \textbf{33} (1984), 400-403.

\bibitem{HS96} T. Hasunuma and Y. Shibata, Isomorphic decomposition and
arc-disjoint spanning trees of Kautz digraphs, IPSJ SIG Notes, 96-AL-51
(1996), 63-70.

\bibitem{KFS01} H. Kawai, N. Fujikake and Y. Shibata, Factorization of de
Bruijn digraphs by cycle-rooted trees, \emph{Inform. Process. Lett.} \textbf{%
77} (2001), \emph{no. 5-6}, 269--275.

\bibitem{H97} M.-C. Heydemann, Cayley graphs and interconnection networks. 
\emph{Graph symmetry} (Montreal, PQ, 1996), 167--224, \emph{NATO Adv. Sci.
Inst. Ser. C Math. Phys. Sci.}, \textbf{497}, Kluwer Acad. Publ., Dordrecht,
1997.
\end{thebibliography}
\end{document}